\documentclass[a4paper,12pt,leqno]{article}
\title{Jumping deformations of complete toric varieties}

\date{}

\author{{\sc Hiroshi Sato}\thanks{Partly supported by the Grant-in-Aid for JSPS Fellows, The Ministry of Education, Science, Sports and Culture, Japan.
\newline
\hspace*{1.5em} {\em $2000$ Mathematics Subject Classification\/}.
Primary 14M25;
Secondary 14J45, 32G05.}}

\newtheorem{Thm}{Theorem}[section]
\newtheorem{Prop}[Thm]{Proposition}

\newtheorem{Def}[Thm]{Definition}
\newtheorem{Rem}[Thm]{Remark}
\newtheorem{Ex}[Thm]{Example}

\newcommand{\proof}{Proof. \quad}

\newcommand{\Hom}{\mathop{\rm Hom}\nolimits}

\newcommand{\G}{\mathop{\rm G}\nolimits}

\begin{document}
\maketitle

\begin{abstract}                                                  

We construct one-parameter complex analytic families whose special fibers are complete toric varieties. Under some assumptions, the general fibers of these families also become toric varieties and we can explicitly describe the corresponding fans from the data of the fans corresponding to the special fibers. Using these families, we give a deformation family for a certain toric weakened Fano $3$-fold. Moreover, we get some examples of toric weakened Fano $4$-folds.

\end{abstract}

\section{Introduction}\label{intro}

\thispagestyle{empty}

\hspace{5mm} It is well-known that the Hirzebruch surface $F_a\ (a\geq 0)$ of degree $a$ is deformed in a one-parameter family to $F_{a-2k}$, where $k$ is a positive integer such that $a-2k\geq 0$. Especially, if $a\equiv a'\ ($mod 2$)$, then $F_a$ and $F_{a'}$ are homeomorphic. In this paper, we generalize this classical result to certain nonsingular complete toric varieties. Namely, for a nonsingular complete toric $d$-fold $V$ which have a toric fibration onto ${\bf P}^1$ such that its general fiber is isomorphic to either ${\bf P}^{d-1}$ or a toric bundle over ${\bf P}^1$, we construct a complex analytic family $\{V_t\}_{t\in {\bf C}}$ such that $V_0\cong V$ and that $\{V_t\}_{t\neq 0}$ are mutually isomorphic. Moreover, under some assumptions, the general fiber of this family can be explicitly described by the data of the fan corresponding to $V$.

As an application of this construction of families, we give deformation families for some toric weakened Fano varieties, that is, nonsingular toric weak Fano varieties which are not Fano but are deformed to Fano varieties. Toric weakened Fano $d$-folds are classified for $d\leq 3$ (see Sato \cite{sato2}). Moreover, we get some examples of toric weakened Fano $4$-folds.

The content of this paper is as follows: In Section \ref{homo}, we review the homogeneous coordinate of a toric varietiy, which is a key to our main result. In Section \ref{main}, we construct complex analytic families of nonsingular complete varieties over {\bf C} as stated above. In Section \ref{scroll}, as an application of the construction, we study deformations among ${\bf P}^{d-1}$-bundles over ${\bf P}^1$. In Section \ref{weakened}, we give some examples of toric weakened Fano $3$-folds and $4$-folds using the families constructed in Section \ref{main}.

The author wishes to thank Professors Tadao Oda, Masa-Nori Ishida, Tadashi Ashikaga, Takeshi Kajiwara and Yasuhiro Nakagawa for advice and encouragement.

\section{Homogeneous coordinates of toric varieties}\label{homo}

\hspace{5mm} In this section, we recall homogeneous coordinates of toric varieties (see Cox \cite{cox1} and Oda \cite{oda2}).

Let $N={\bf Z}^{d}$ with elements regarded as column vectors, $M:=\Hom_{\bf Z}(N,{\bf Z})$, $N_{\bf R}:=N\otimes {\bf R}$, $M_{\bf R}:=M\otimes {\bf R}$ and $\Sigma$ a fan in $N$. Throughout this paper, we mean by a cone a nonsingular rational cone and by a fan in $N$ a nonsingular fan which contains at least one $d$-dimensional cone. For $0\leq i\leq d$, we put $\Sigma(i):=\left\{ \sigma\in\Sigma\; |\;\dim\sigma=i\right\}$. Each $\tau\in\Sigma(1)$ determines a unique element $e(\tau)\in N$ which generates the semigroup $\tau\cap N$. We call
$$\G(\Sigma):=\left\{e(\tau)\in N\; |\;\tau\in\Sigma(1)\right\}$$
the set of primitive generators of $\Sigma$. Put $\G(\sigma):=\sigma\cap\G(\Sigma)$. We introduce variables $\{\mathcal{X}_{\rho}\; |\; \rho\in\G(\Sigma)\}$ and consider the polynomial ring $S:={\bf C}\left[ \mathcal{X}_{\rho}\; |\; \rho\in\G(\Sigma)\right]$, which we call the {\em homogeneous coordinate ring} of the nonsingular toric $d$-fold $V$ corresponding to $\Sigma$. Put
$$Z:=\left\{ X=(X_{\rho})_{\rho\in\G(\Sigma)}\in{\bf C}^{\G(\Sigma)}\; \left|\; \prod_{\rho\in\G(\Sigma)\setminus\G(\sigma)}X_{\rho}=0\mbox{ for any }\sigma\in\Sigma\right.\right\}\subset {\bf C}^{\G(\Sigma)}.$$
On the other hand, by the exact sequence
$$0\rightarrow M\rightarrow {\bf Z}^{\G(\Sigma)}\rightarrow {\rm Pic}(V)\rightarrow 0,$$
we have an exact sequence
$$1\rightarrow G:={\rm Hom}_{\bf Z}\left( {\rm Pic}(V),{\bf C}^{\times}\right) \rightarrow \left({\bf C}^{\times}\right)^{\G(\Sigma)}\rightarrow T_N\rightarrow 1.$$
Since $\left({\bf C}^{\times}\right)^{\G(\Sigma)}$ acts naturally on ${\bf C}^{\G(\Sigma)}$, the subgroup $G\subset\left({\bf C}^{\times}\right)^{\G(\Sigma)}$ acts on ${\bf C}^{\G(\Sigma)}$ as
$$gt:=\left( g\left( \left[D_{\rho}\right]\right)t_{\rho}\right)_{\rho\in\G(\Sigma)},$$
where $g\in G$, $t=\left( t_{\rho}\right)_{\rho\in\G(\Sigma)}\in{\bf C}^{\G(\Sigma)}$ and $[D_{\rho}]\in{\rm Pic}(V)$ is the class of the $T_N$-invariant prime divisor $D_{\rho}$ corresponding to $\rho$. In this setting, the following holds.

\begin{Prop}[Cox \cite{cox1}, Theorem 2.1]$\!\!\!${\bf .}\label{homoquot}
The subset ${\bf C}^{\G(\Sigma)}\setminus Z\subset{\bf C}^{\G(\Sigma)}$ is invariant under the action of $G$, and $V$ is the geometric quotient of ${\bf C}^{\G(\Sigma)}\setminus Z$ by $G$. We denote ${\bf C}^{\G(\Sigma)}\setminus Z$ by ${\bf U}(\Sigma)$.

\end{Prop}

We need the following proposition for this description.

\begin{Prop}[Cox \cite{cox1}, Theorem 2.1]$\!\!\!${\bf .}
For any $\sigma \in\Sigma$, we have
$$U_{\sigma}\cong\left({\bf U}(\Sigma)_{\sigma}:=\left.\left\{ X=\left( X_{\rho}\right)_{\rho\in\G(\Sigma)}\in{\bf U}(\Sigma)\; \left|\; \prod_{\rho\in\G(\Sigma)\setminus\G(\sigma)}X_{\rho}\neq 0\right.\right\}\right)\right/G,$$
where $U_{\sigma}\subset V$ is the affine toric subvariety corresponding to $\sigma$.

\end{Prop}

\section{Constructions of families}\label{main}

\hspace{5mm} In this section, we construct one-parameter complex analytic families whose fibers are nonsingular complete varieties. Especially, the special fibers are nonsingular complete toric varieties. This is a generalization of the classical results on deformations among Hirzebruch surfaces.


Let $\widetilde{N}:=\left\{ {\bf n} \in N \,|\, \mbox{the }d\mbox{-th coordinate of } {\bf n} \mbox{ is } 0\right\}$ and $\widetilde{\Sigma}$ a complete fan in $\widetilde{N}$. For a complete fan $\Sigma$ in $N$ containing $\widetilde{\Sigma}$ as a subfan, we define subfans of $\Sigma$ as follows:
$$\Sigma^{+}:=\left\{ \sigma \in \Sigma \,|\, \mbox{the }d\mbox{-th coordinate of } {\bf n} \mbox{ is nonnegative for any } {\bf n} \in \sigma \right\}$$
$$\Sigma^{-}:=\left\{ \sigma \in \Sigma \,|\, \mbox{the }d\mbox{-th coordinate of } {\bf n} \mbox{ is nonpositive for any } {\bf n} \in \sigma \right\}$$
Then we have $\widetilde{\Sigma}=\Sigma^{+}\cap\Sigma^{-}$. We denote by $V$ (resp. $V^{+}$, $V^{-}$, $\widetilde{V}$) the nonsingular toric variety corresponding to the fan $\Sigma$ (resp. $\Sigma^{+}$, $\Sigma^{-}$, $\widetilde{\Sigma}$).

\begin{Rem}

{\rm
$V$ has a toric fibration $V\rightarrow{\bf P}^1$ whose general fiber is isomorphic to $\widetilde{V}$.
}

\end{Rem}

In the above situation, let
$$\G(\widetilde{\Sigma})=\left\{ {\bf e}_1,\ldots,{\bf e}_{d-1},{\bf a}_{1},\ldots,{\bf a}_{\rho}\right\},\ \G(\Sigma^{+})=\left\{ {\bf b}_{1},\ldots,{\bf b}_{m}\right\}\cup\G(\widetilde{\Sigma}),$$
$$\G(\Sigma^{-})=\left\{ {\bf c}_{1},\ldots,{\bf c}_{n}\right\}\cup\G(\widetilde{\Sigma}),$$
$\left\{ {\bf e}_1,\ldots,{\bf e}_{d-1},{\bf b}_1\right\}$ the standard basis for $N$ and
\[
\left( {\bf a}_1,\ldots,{\bf a}_{\rho},{\bf b}_2,\ldots,{\bf b}_{m},{\bf c}_1,\ldots,{\bf c}_{n}\right) =
\]
\[
\pmatrix{
a_{1,1} & \cdots & a_{\rho,1} & b_{2,1} & \cdots & b_{m,1} & c_{1,1} & \cdots & c_{n,1} \cr
\vdots & \ddots & \vdots & \vdots & \ddots & \vdots & \vdots & \ddots & \vdots \cr
a_{1,d} & \cdots & a_{\rho,d} & b_{2,d} & \cdots & b_{m,d} & c_{1,d} & \cdots & c_{n,d}}.
\]
Suppose that $\widetilde{V}$ is isomorphic to either ${\bf P}^{d-1}$ or a toric bundle over ${\bf P}^1$. If $\widetilde{V}$ is isomorphic to a toric bundle over ${\bf P}^1$, suppose that the $T_N$-invariant prime divisors on $\widetilde{V}$ corresponding to ${\bf e}_1$ and ${\bf a}_1$ correspond to fibers. Suppose further that $\left\{ {\bf e}_1,\ldots,{\bf e}_{d-1},{\bf b}_1\right\}$ generates a $d$-dimensional cone in $\Sigma^{+}$, while $\left\{ {\bf e}_1,\ldots,{\bf e}_{d-1},{\bf c}_1\right\}$ generates a $d$-dimensional cone in $\Sigma^{-}$. For a nonnegative integer $k$, we construct a complex analytic family.

Since $\left\{ {\bf a}_{1},\ldots,{\bf a}_{\rho}\right\}\subset\G(\widetilde{\Sigma})$, we have $a_{1,d}=\cdots =a_{\rho,d}=0$. We have $c_{1,d}=-1$, by the assumption that $\left\{ {\bf e}_1,\ldots,{\bf e}_{d-1},{\bf c}_1\right\}$ generates a $d$-dimensional cone in $\Sigma^{-}$.

Let $D_1,\ldots,D_{d-1},A_1,\ldots,A_{\rho},B_1,\ldots,B_{m},C_1,\ldots,C_{n}$ be the $T_N$-invariant prime divisors corresponding to ${\bf e}_1,\ldots,{\bf e}_{d-1},{\bf a}_1,\ldots,{\bf a}_{\rho},{\bf b}_1,\ldots,{\bf b}_{m},{\bf c}_1,\ldots,{\bf c}_n$, respectively. Then by computing the divisors of the rational functions ${\bf e}({\bf e}^{\ast}_1),\ldots,{\bf e}({\bf e}^{\ast}_{d-1}),{\bf e}({\bf b}^{\ast}_1)\in {\bf C}(V)$, where $\{ {\bf e}^{\ast}_1,\ldots,{\bf e}^{\ast}_{d-1},{\bf b}^{\ast}_1\}\subset M$ is the dual basis of $\{ {\bf e}_1,\ldots,{\bf e}_{d-1},{\bf b}_1\}$, we have
$$\ D_1+a_{1,1}A_1+\cdots+a_{\rho,1}A_{\rho}+b_{2,1}B_2+\cdots+b_{m,1}B_m+c_{1,1}C_1+\cdots+c_{n,1}C_n=0,$$
$$D_2+a_{1,2}A_1+\cdots+a_{\rho,2}A_{\rho}+b_{2,2}B_2+\cdots+b_{m,2}B_m+c_{1,2}C_1+\cdots+c_{n,2}C_n=0,$$
$$\vdots$$
$$D_{d-1}+a_{1,d-1}A_1+\cdots+a_{\rho,d-1}A_{\rho}+b_{2,d-1}B_2+\cdots+b_{m,d-1}B_m+c_{1,d-1}C_1+\cdots+c_{n,d-1}C_n=0,$$
$$B_1+b_{2,d}B_2+\cdots+b_{m,d}B_m-C_1+c_{2,d}C_2+\cdots+c_{n,d}C_n=0$$
in ${\rm Pic}(V)$, respectively. Using these equalities, we calculate the homogeneous coordinates of $V$, $V^{+}$, $V^{-}$ and $\widetilde{V}$.

Let $\left( X_1,\ldots,X_{d-1},Y_1,\ldots,Y_{\rho},Z_{1},\ldots,Z_{m},W_{1},\ldots,W_{n}\right) \in {\bf U}(\Sigma)$ be the homogeneous coordinate of $V$ corresponding to ${\bf e}_1,\ldots,{\bf e}_{d-1},{\bf a}_1,\ldots,{\bf a}_{\rho},{\bf b}_1,\ldots,{\bf b}_{m},{\bf c}_1,\ldots,{\bf c}_n$, respectively. Then the action of $G:={\rm Hom}_{\bf Z}({\rm Pic}(V),{\bf C}^{\times})$ on ${\bf U}(\Sigma)$ is as follows: $g\in G$ acts as
\begin{equation}
\left( g\left( -\left( a_{1,1}A_1+\cdots+a_{\rho,1}A_{\rho}+b_{2,1}B_2+\cdots+b_{m,1}B_m+c_{1,1}C_1+\cdots+c_{n,1}C_n\right) \right)X_1,\right.
\end{equation}
$$\ldots,g\left( -\left( a_{1,d-1}A_1+\cdots+a_{\rho,d-1}A_{\rho}+b_{2,d-1}B_2+\cdots+b_{m,d-1}B_m+\right.\right.$$
$$\left.\left.c_{1,d-1}C_1+\cdots+c_{n,d-1}C_n\right) \right)X_{d-1},$$
$$g(A_1)Y_1,\ldots,g(A_{\rho})Y_{\rho},$$
$$g\left( -\left( b_{2,d}B_2+\cdots+b_{m,d}B_m-C_1+c_{2,d}C_2+\cdots+c_{n,d}C_n\right) \right) Z_{1},$$
$$\left. g(B_2)Z_{2},\ldots,g(B_m)Z_{m},g(C_1)W_{1},\cdots,g(C_n)W_{n}\right).$$
Let $\varphi:(X_1,\ldots,X_{d-1},Y_1,\ldots,Y_{\rho},Z_{1},\ldots,Z_{m},W_{1},\ldots,W_{n})\mapsto(X^{+}_1,\ldots,X^{+}_{d-1},Y^{+}_1,\ldots,Y^{+}_{\rho},$ $Z^{+}_{1},\ldots,Z^{+}_{m})$ be a surjective morphism from $\bigcup_{\sigma\in\Sigma^+}{\bf U}(\Sigma)_{\sigma}\subset{\bf U}(\Sigma)$ to ${\bf U}(\Sigma^{+})$ given by
\begin{equation}
\left( X^{+}_1,\ldots,X^{+}_{d-1},Y^{+}_1,\ldots,Y^{+}_{\rho},Z^{+}_{1},\ldots,Z^{+}_{m}\right)=
\end{equation}
$$\left( X_1W^{c_{1,1}}_1\cdots W^{c_{n,1}}_n,\ldots,X_{d-1}W^{c_{1,d-1}}_1\cdots W^{c_{n,d-1}}_n,\right.$$
$$\left. Y_1,\ldots,Y_{\rho},Z_{1}W^{c_{1,d}}_1\cdots W^{c_{n,d}}_n,Z_2,\ldots,Z_{m}\right),$$
where $\left( X^{+}_1,\ldots,X^{+}_{d-1},Y^{+}_1,\ldots,Y^{+}_{\rho},Z^{+}_{1},\ldots,Z^{+}_{m}\right)$ is the homogeneous coordinate of $V^{+}\cong{\bf U}(\Sigma^+)/G^+$ with $G^{+}:={\rm Hom}_{\bf Z}({\rm Pic}(V^{+}),{\bf C}^{\times})$ corresponding to ${\bf e}_1,\ldots,{\bf e}_{d-1},{\bf a}_1,\ldots,{\bf a}_{\rho},{\bf b}_1,$ $\ldots,{\bf b}_{m}$, respectively. $\varphi$ is well-defined, since $W_1,\ldots,W_n\neq 0$ on $\bigcup_{\sigma\in\Sigma^+}{\bf U}(\Sigma)_{\sigma}$. Moreover, since $\varphi$ is compartible with the action of $G$ and $G^{+}$ by $(1)$, $\varphi$ induces the isomorphism $\widetilde{\varphi}:\left(\bigcup_{\sigma\in\Sigma^{+}}{\bf U}(\Sigma)_{\sigma}\right)/G\subset V\rightarrow V^+$. Similarly, the morphism $\psi:(X_1,\ldots,X_{d-1},Y_1,\ldots,Y_{\rho},Z_{1},$ $\ldots,Z_{m},W_{1},\ldots,W_{n})\mapsto(X^{-}_1,\ldots,X^{-}_{d-1},Y^{-}_1,\ldots,Y^{-}_{\rho},W^{-}_{1},\ldots,W^{-}_{n})$ from \ \ $\bigcup_{\sigma\in\Sigma^{-}}{\bf U}(\Sigma)_{\sigma}\subset{\bf U}(\Sigma)$ to ${\bf U}(\Sigma^{-})$ given by
\begin{equation}
\left( X^{-}_1,\ldots,X^{-}_{d-1},Y^{-}_1,\ldots,Y^{-}_{\rho},W^{-}_{1},\ldots,W^{-}_{n}\right)=
\end{equation}
$$\left( X_1Z^{c_{1,1}}_1Z^{b_{2,d}c_{1,1}+b_{2,1}}_2\cdots Z^{b_{m,d}c_{1,1}+b_{m,1}}_{m},\ldots,\right.$$
$$\left.X_{d-1}Z^{c_{1,d-1}}_1Z^{b_{2,d}c_{1,d-1}+b_{2,d-1}}_2\cdots Z^{b_{m,d}c_{1,d-1}+b_{m,d-1}}_{m},Y_1,\ldots,Y_{\rho},\right.$$
$$\left.Z^{-1}_1Z^{-b_{2,d}}_2\cdots Z^{-b_{m,d}}_mW_1,W_2,\ldots,W_{n}\right)$$
induces the isomorphism $\widetilde{\psi}:\left(\bigcup_{\sigma\in\Sigma^{-}}{\bf U}(\Sigma)_{\sigma}\right)/G\subset V\rightarrow V^-\cong{\bf U}(\Sigma^{-})/G^-$, where $G^{-}:={\rm Hom}_{\bf Z}({\rm Pic}(V^{-}),{\bf C}^{\times})$ and $\left( X^{-}_1,\ldots,X^{-}_{d-1},Y^{-}_1,\ldots,Y^{-}_{\rho},W^{-}_{1},\ldots,W^{-}_{n}\right)$ is the homogeneous coordinate of $V^{-}$ corresponding to ${\bf e}_1,\ldots,{\bf e}_{d-1},{\bf a}_1,\ldots,{\bf a}_{\rho},{\bf c}_1,\ldots,{\bf c}_{n}$, respectively.

$g^{+}\in G^{+}$ and $g^{-}\in G^{-}$ act on ${\bf U}(\Sigma^{+})$ and ${\bf U}(\Sigma^{-})$ as
\begin{equation}
\left( g^{+}\left( -\left( a_{1,1}A_1+\cdots+a_{\rho,1}A_{\rho}+b_{2,1}B_2+\cdots+b_{m,1}B_m\right) \right)X^{+}_1,\ldots,\right.
\end{equation}
$$g^{+}\left( -\left( a_{1,d-1}A_1+\cdots+a_{\rho,d-1}A_{\rho}+b_{2,d-1}B_2+\cdots+b_{m,d-1}B_m\right) \right)X^{+}_{d-1},$$
$$g^{+}(A_1)Y^{+}_1,\ldots,g^{+}(A_{\rho})Y^{+}_{\rho},g^{+}\left( -\left( b_{2,d}B_2+\cdots+b_{m,d}B_m\right) \right) Z^{+}_{1},$$
$$\left. g^{+}(B_2)Z^{+}_{2},\ldots,g^{+}(B_m)Z^{+}_{m}\right) $$
$$\mbox{and}$$
$$\left( g^{-}\left( -\left( a_{1,1}A_1+\cdots+a_{\rho,1}A_{\rho}+(c_{2,1}+c_{1,1}c_{2,d})C_2+\cdots+(c_{n,1}+c_{1,1}c_{n,d})C_n\right) \right)X^{-}_1,\ldots,\right.$$
$$g^{-}\left( -\left( a_{1,d-1}A_1+\cdots+a_{\rho,d-1}A_{\rho}+(c_{2,d-1}+c_{1,d-1}c_{2,d})C_2+\cdots+\right.\right.$$
$$\left.\left.(c_{n,d-1}+c_{1,d-1}c_{n,d})C_n\right) \right)X^{-}_{d-1},g^{-}(A_1)Y^{-}_1,\ldots,g^{-}(A_{\rho})Y^{-}_{\rho},$$
$$\left. g^{-}\left( -(c_{2,d}C_2+\cdots+c_{n,d}C_n)\right) W^{-}_{1},g^{-}(C_2)W^{-}_{2},\cdots,g^{-}(C_n)W^{-}_{n}\right),$$
respectively. So, similarly as $(2)$ and $(3)$, we have isomorphisms $(X^{+}_1,\ldots,X^{+}_{d-1},Y^{+}_1,\ldots,$ $Y^{+}_{\rho},Z^{+}_{1},\ldots,Z^{+}_{m})\mapsto(x^{+}_1,$ $\ldots,x^{+}_{d-1},y^{+}_1,\ldots,y^{+}_{\rho},z)$ from $\left(\bigcup_{\sigma\in\widetilde{\Sigma}}{\bf U}(\Sigma^+)_{\sigma}\right)/G^{+}\subset V^{+}$ to $\widetilde{V}\times {\bf C}^{\times}$ and $(X^{-}_1,\ldots,X^{-}_{d-1},Y^{-}_1,\ldots,Y^{-}_{\rho},$ $W^{-}_{1},\ldots,W^{-}_{n})\mapsto(x^{-}_1,\ldots,x^{-}_{d-1},y^{-}_1,\ldots,y^{-}_{\rho},w)$ from $\left(\bigcup_{\sigma\in\widetilde{\Sigma}}{\bf U}(\Sigma^-)_{\sigma}\right)/G^{-}\subset V^{-}$ to $\widetilde{V}\times {\bf C}^{\times}$ given by
\begin{equation}
\left( x^{+}_1,\ldots,x^{+}_{d-1},y^{+}_1,\ldots,y^{+}_{\rho},z\right)=
\end{equation}
$$\left( X^{+}_1(Z^{+}_2)^{b_{2,1}}\cdots (Z^{+}_m)^{b_{m,1}},\ldots,X^{+}_{d-1}(Z^{+}_2)^{b_{2,d-1}}\cdots (Z^{+}_m)^{b_{m,d-1}},\right.$$
$$\left. Y^{+}_1,\ldots,Y^{+}_{\rho},Z^{+}_{1}(Z^{+}_2)^{b_{2,d}}\cdots (Z^{+}_m)^{b_{m,d}}\right)$$
and
\begin{equation}
\left( x^{-}_1,\ldots,x^{-}_{d-1},y^{-}_1,\ldots,y^{-}_{\rho},w\right)=
\end{equation}
$$\left( X^{-}_1(W^{-}_2)^{c_{2,1}+c_{1,1}c_{2,d}}\cdots (W^{-}_n)^{c_{n,1}+c_{1,1}c_{n,d}},\ldots,\right.$$
$$\left.X^{-}_{d-1}(W^{-}_2)^{c_{2,d-1}+c_{1,d-1}c_{2,d}}\cdots(W^{-}_n)^{c_{n,d-1}+c_{1,d-1}c_{n,d}},Y^{-}_1,\ldots,Y^{-}_{\rho},\right.$$
$$\left.W^{-}_{1}(W^{-}_2)^{-c_{2,d}}\cdots (W^{-}_n)^{-c_{n,d}}\right),$$
respectively, where $( x^{+}_1,\ldots,x^{+}_{d-1},y^{+}_1,\ldots,y^{+}_{\rho})$ and $( x^{-}_1,\ldots,x^{-}_{d-1},y^{-}_1,\ldots,y^{-}_{\rho})$ are homogeneous coordinates of $\widetilde{V}$, while $z,w\in{\bf C}^{\times}$. These two coordinates are related as follows:
$$x^{+}_1=x^{-}_1w^{c_{1,1}},\ldots,x^{+}_{d-1}=x^{-}_{d-1}w^{c_{1,d-1}},y^{+}_1=y^{-}_1,\ldots,y^{+}_{\rho}=y^{-}_{\rho},z=1/w.$$
We construct a one-parameter family of toric varieties parameterized by $t\in{\bf C}$ by changing this relation: Let $\{ V_t\}_{t\in{\bf C}}$ be the family we obtain by patching $V^{+}$ and $V^{-}$ along $\widetilde{V}$ by the automorphism $(x^{-}_1,\ldots,x^{-}_{d-1},y^{-}_1,\ldots,y^{-}_{\rho},w)\mapsto(x^{+}_1,\ldots,x^{+}_{d-1},y^{+}_1,\ldots,y^{+}_{\rho},z)$ defined by
\begin{equation}
x^{+}_1=x^{-}_1w^{c_{1,1}}+ty^{-}_1w^{k},x^{+}_{2}=x^{-}_{2}w^{c_{1,2}},\ldots,x^{+}_{d-1}=x^{-}_{d-1}w^{c_{1,d-1}},
\end{equation}
$$y^{+}_1=y^{-}_1,\ldots,y^{+}_{\rho}=y^{-}_{\rho},z=1/w.$$
This is well-defined, since $D_1=A_1$ in ${\rm Pic}(\widetilde{V})$ and the combinatorial structures of the neighborhoods of ${\bf e}_1$ and ${\bf a}_1$ in $\widetilde{\Sigma}$ are equivalent by the assumption $\widetilde{V}$ is isomorphic to either ${\bf P}^{d-1}$ or a toric bundle over ${\bf P}^1$. Thus, we have the following.

\begin{Thm}\label{maindeform}

$\{ V_t\}_{t\in{\bf C}}$ is a complex analytic family whose special fiber $V_0$ is isomorphic to $V$.

\end{Thm}


Next, we calculate the general fibers of this family under some assumptions. We introduce some notation.

For any ${\bf q}=\left( q_1,\ldots,q_{d-1}\right) \in {\bf Z}^{d-1}$ we can define a complete fan ${\bf q}^{-}\Sigma$ in $N$ as follows:
$${\bf q}^{-}\Sigma:=\Sigma^{+}\cup\left\{ {\bf q}^{-}\sigma \;|\; \sigma\in\Sigma^{-} \right\},$$
where ${\bf q}^{-}\sigma$ is the image of $\sigma$ under the automorphism of $N_{\bf R}$ corresponding to the matrix acting from the left on the elements of $N={\bf Z}^{d}$ regarded as column vectors
\[
\pmatrix{
1 & 0 & \cdots & 0 & q_1 \cr
0 & 1 & \cdots & 0 & q_2 \cr
\vdots & \vdots & \ddots & \vdots & \vdots \cr
0 & 0 & \cdots & 1 & q_{d-1} \cr
0 & 0 & \cdots & 0 & 1}.
\]
We denotes by $q^{-}V$ the nonsingular toric $d$-fold corresponding to the fan $q^{-}\Sigma$.

\begin{Thm}\label{maindeform2}

For any $t$ in ${\bf C}^{\times}:={\bf C}\setminus\{0\}$, we have
$$V_t\cong\left(2k,-ka_{1,2},\ldots,-ka_{1,d-1}\right)^{-}V,$$
if the following conditions are satisfied$:$
\begin{enumerate}
\item $b_{2,1}=\cdots=b_{m,1}=0$ and
\item $kc_{1,d}+c_{1,1}\geq 0,\ldots,kc_{n,d}+c_{n,1}\geq 0.$
\end{enumerate}

\end{Thm}

\proof
Let $t\neq 0$.

We can define an automorphism $(x^{+}_1,\ldots,x^{+}_{d-1},y^{+}_1,\ldots,y^{+}_{\rho},z)\mapsto(\hat{x}^{+}_1,\ldots,\hat{x}^{+}_{d-1},\hat{y}^{+}_1,$ $\ldots,\hat{y}^{+}_{\rho},\hat{z})$ of $\widetilde{V}\times{\bf C}^{\times}$ by
\begin{equation}
\hat{x}^{+}_1:=x^{+}_1z^{k}-ty^{+}_1,\hat{x}^{+}_2:=x^{+}_2,\ldots,\hat{x}^{+}_{d-1}:=x^{+}_{d-1},
\end{equation}
$$\hat{y}^{+}_1:=tx^{+}_1,\hat{y}^{+}_2:=y^{+}_2,\ldots,\hat{y}^{+}_{\rho}:=y^{+}_{\rho},\hat{z}:=z.$$
In fact, since $\widetilde{V}$ is isomorphic to either ${\bf P}^{d-1}$ or a toric bundle over ${\bf P}^1$, this is a morphism, and we can easily construct the inverse of this morphism. By the automorphism $(x^{-}_1,\ldots,x^{-}_{d-1},y^{-}_1,\ldots,y^{-}_{\rho},w)\mapsto(x^{+}_1,\ldots,x^{+}_{d-1},y^{+}_1,\ldots,y^{+}_{\rho},z)$ given by the relations $(7)$, we have
$$\hat{x}^{+}_1=(x^{-}_1w^{c_{1,1}}+ty^{-}_1w^{k})w^{-k}-ty^{-}_1=x^{-}_1w^{c_{1,1}-k},$$
$$\hat{x}^{+}_2=x^{-}_{2}w^{c_{1,2}},\ldots,\hat{x}^{+}_{d-1}=x^{-}_{d-1}w^{c_{1,d-1}},$$
$$\hat{y}^{+}_1=t(x^{-}_1w^{c_{1,1}}+ty^{-}_1w^{k}),\ \hat{y}^{+}_2=y^{-}_2,\ldots,\hat{y}^{+}_{\rho}=y^{-}_{\rho},\ \hat{z}=1/w.$$
By considering the action of $\widetilde{G}$ on $x^{+}_{1},\ldots,x^{+}_{d-1},y^{+}_1$, these relations are equivalent to
\begin{equation}
\hat{x}^{+}_1=x^{-}_1w^{c_{1,1}-2k},\ \hat{x}^{+}_2=x^{-}_{2}w^{c_{1,2}+ka_{1,2}},\ldots,\hat{x}^{+}_{d-1}=x^{-}_{d-1}w^{c_{1,d-1}+ka_{1,d-1}},
\end{equation}
$$\hat{y}^{+}_1=t(x^{-}_1w^{c_{1,1}-k}+ty^{-}_1),\ \hat{y}^{+}_2=y^{-}_2,\ldots,\hat{y}^{+}_{\rho}=y^{-}_{\rho},\ \hat{z}=1/w.$$
Let
\begin{equation}
\hat{x}^{-}_1:=x^{-}_1,\ldots,\hat{x}^{-}_{d-1}:=x^{-}_{d-1},
\end{equation}
$$\hat{y}^{-}_1:=t(x^{-}_1w^{c_{1,1}-k}+ty^{-}_1),\ \hat{y}^{-}_2:=y^{-}_2,\ldots,\hat{y}^{-}_{\rho}:=y^{-}_{\rho},\ \hat{w}:=w.$$
Then $(x^{-}_1,\ldots,x^{-}_{d-1},y^{-}_1,\ldots,y^{-}_{\rho},w)\mapsto(\hat{x}^{-}_1,\ldots,\hat{x}^{-}_{d-1},\hat{y}^{-}_1,\ldots,\hat{y}^{-}_{\rho},\hat{w})$ determines an automorphism of $\widetilde{V}\times{\bf C}^{\times}$, and in terms of this new coordinate, the automorphism $(x^{-}_1,\ldots,x^{-}_{d-1},$ $y^{-}_1,\ldots,y^{-}_{\rho},w)\mapsto(\hat{x}^{+}_1,\ldots,\hat{x}^{+}_{d-1},\hat{y}^{+}_1,\ldots,\hat{y}^{+}_{\rho},\hat{z})$ given by the relations $(9)$ is described as the automorphism $(\hat{x}^{-}_1,\ldots,\hat{x}^{-}_{d-1},\hat{y}^{-}_1,\ldots,\hat{y}^{-}_{\rho},\hat{w})\mapsto(\hat{x}^{+}_1,\ldots,\hat{x}^{+}_{d-1},\hat{y}^{+}_1,\ldots,\hat{y}^{+}_{\rho},\hat{z})$ given by
\begin{equation}
\hat{x}^{+}_1={\hat{x}_1}^{-}\hat{w}^{c_{1,1}-2k},\ \hat{x}^{+}_2={\hat{x}_2}^{-}\hat{w}^{c_{1,2}+ka_{1,2}},\ldots,\hat{x}^{+}_{d-1}={\hat{x}_{d-1}}^{-}\hat{w}^{c_{1,d-1}+ka_{1,d-1}},
\end{equation}
$$\hat{y}^{+}_1=\hat{y}^{-}_1,\ldots,\hat{y}^{+}_{\rho}={\hat{y}_{\rho}}^{-},\ \hat{z}=1/\hat{w}.$$
We can show that the automorphisms $(x^{+}_1,\ldots,x^{+}_{d-1},y^{+}_1,\ldots,y^{+}_{\rho},z)\mapsto(\hat{x}^{+}_1,\ldots,\hat{x}^{+}_{d-1},\hat{y}^{+}_1,$ $\ldots,\hat{y}^{+}_{\rho},\hat{z})$ in $(8)$ and $(x^{-}_1,\ldots,x^{-}_{d-1},y^{-}_1,\ldots,y^{-}_{\rho},w)\mapsto(\hat{x}^{-}_1,\ldots,\hat{x}^{-}_{d-1},\hat{y}^{-}_1,\ldots,\hat{y}^{-}_{\rho},\hat{w})$ in $(10)$ of $\widetilde{V}\times{\bf C}^{\times}$ are extended to automorphisms of $V^{+}$ and $V^{-}$, respectively as follows: Put
$$\hat{X}^{+}_1:=X^{+}_1(Z^{+}_1)^{k}(Z^{+}_2)^{kb_{2,d}-b_{2,1}}\cdots (Z^{+}_m)^{kb_{m,d}-b_{m,1}}-tY^{+}_1(Z^{+}_2)^{-b_{2,1}}\cdots (Z^{+}_m)^{-b_{m,1}},$$
$$\hat{X}^{+}_2:=X^{+}_2,\ldots,\hat{X}^{+}_{d-1}:=X^{+}_{d-1},\ \hat{Y}^{+}_1:=tX^{+}_1(Z^{+}_2)^{b_{2,1}}\cdots (Z^{+}_m)^{b_{m,1}},$$
$$\hat{Y}^{+}_2:=Y^{+}_2,\ldots,\hat{Y}^{+}_{\rho}:=Y^{+}_{\rho},\ \hat{Z}^{+}_1:=Z^{+}_1,\ldots,\hat{Z}^{+}_{m}:=Z^{+}_{m}.$$
By the assumption $b_{2,1}=\cdots=b_{m,1}=0$, this defines an automorphism $(X^{+}_1,\ldots,X^{+}_{d-1},$ $Y^{+}_1,\ldots,Y^{+}_{\rho},Z^{+}_{1},\ldots,Z^{+}_{m})\mapsto(\hat{X}^{+}_1,\ldots,\hat{X}^{+}_{d-1},\hat{Y}^{+}_1,\ldots,\hat{Y}^{+}_{\rho},\hat{Z}^{+}_{1},\ldots,\hat{Z}^{+}_{m})$ of $V^{+}$, and obviously the restriction of this automorphism through the isomorphisms $(X^{+}_1,\ldots,X^{+}_{d-1},Y^{+}_1,$ $\ldots,Y^{+}_{\rho},Z^{+}_{1},\ldots,Z^{+}_{m})\mapsto(x^{+}_1,\ldots,x^{+}_{d-1},y^{+}_1,\ldots,y^{+}_{\rho},z)$ and $(\hat{X}^{+}_1,\ldots,\hat{X}^{+}_{d-1},\hat{Y}^{+}_1,\ldots,\hat{Y}^{+}_{\rho},\hat{Z}^{+}_{1},$ $\ldots,\hat{Z}^{+}_{m})\mapsto(\hat{x}^{+}_1,\ldots,\hat{x}^{+}_{d-1},\hat{y}^{+}_1,\ldots,\hat{y}^{+}_{\rho},\hat{z})$ defined by the equalities $(5)$ from ${\bf U}(\widetilde{\Sigma})/G^{+}\subset V^{+}$ to $\widetilde{V}\times{\bf C}^{\times}$ is the automorphism $(x^{+}_1,\ldots,x^{+}_{d-1},y^{+}_1,\ldots,y^{+}_{\rho},z)\mapsto(\hat{x}^{+}_1,\ldots,\hat{x}^{+}_{d-1},\hat{y}^{+}_1,\ldots,\hat{y}^{+}_{\rho},$ $\hat{z})$ corresponding to the equalities $(8)$. Similarly, by the assumption $kc_{1,d}+c_{1,1}\geq 0,\ldots,kc_{n,d}+c_{n,1}\geq 0$, by putting
$$\hat{X}^{-}_1:=X^{-}_1,\ldots,\hat{X}^{-}_{d-1}:=X^{-}_{d-1},$$
$$\hat{Y}^{-}_1:=tX^{-}_1(W^{-}_1)^{kc_{1,d}+c_{1,1}}\cdots (W^{-}_n)^{kc_{n,d}+c_{n,1}},$$
$$\hat{Y}^{-}_2:=Y^{-}_2,\ldots,\hat{Y}^{-}_{\rho}:=Y^{-}_{\rho},\ \hat{W}^{-}_1:=W^{-}_1,\ldots,\hat{W}^{-}_{n}:=W^{-}_{n},$$
we get an automorphism $(X^{-}_1,\ldots,X^{-}_{d-1},Y^{-}_1,\ldots,Y^{-}_{\rho},W^{-}_{1},\ldots,W^{-}_{n})\mapsto(\hat{X}^{-}_1,\ldots,\hat{X}^{-}_{d-1},$ $\hat{Y}^{-}_1,\ldots,\hat{Y}^{-}_{\rho},\hat{W}^{-}_{1},\ldots,\hat{W}^{-}_{n})$ of $V^{-}$ whose restriction through the isomorphisms $(X^{-}_1,\ldots,$ $X^{-}_{d-1},Y^{-}_1,\ldots,Y^{-}_{\rho},W^{-}_{1},\ldots,W^{-}_{n})\mapsto(x^{-}_1,\ldots,x^{-}_{d-1},y^{-}_1,\ldots,y^{-}_{\rho},w)$ and $(\hat{X}^{-}_1,\ldots,\hat{X}^{-}_{d-1},\hat{Y}^{-}_1,$ $\ldots,\hat{Y}^{-}_{\rho},\hat{W}^{-}_{1},\ldots,\hat{W}^{-}_{n})\mapsto(\hat{x}^{-}_1,\ldots,\hat{x}^{-}_{d-1},\hat{y}^{-}_1,\ldots,\hat{y}^{-}_{\rho},\hat{w})$ defined by the equalities $(6)$ from ${\bf U}(\widetilde{\Sigma})/G^{-}\subset V^{-}$ to $\widetilde{V}\times{\bf C}^{\times}$ is the automorphism $(x^{-}_1,\ldots,x^{-}_{d-1},y^{-}_1,\ldots,y^{-}_{\rho},w)\mapsto(\hat{x}^{-}_1,\ldots,$ $\hat{x}^{-}_{d-1},\hat{y}^{-}_1,\ldots,\hat{y}^{-}_{\rho},\hat{w})$ corresponding to the equalities $(10)$.

Next, to show $V_{t}\cong \left(2k,-ka_{1,2},\ldots,-ka_{1,d-1}\right)^{-}V$ for any $t\neq 0$, we have to investigate the action of $G^{-}$ on ${\bf U}(\Sigma^{-})$. However the action $(4)$ is obviously equivalent to the following:
\begin{equation}
\left( g^{-}\left( -\left( a_{1,1}A_1+\cdots+a_{\rho,1}A_{\rho}+(-c_{2,1}-c_{1,1}c_{2,d})C_2+\cdots\right.\right.\right.
\end{equation}
$$\left.\left.\left.+(-c_{n,1}-c_{1,1}c_{n,d})C_n\right) \right)X^{-}_1,\ldots,\right.$$
$$g^{-}\left( -\left( a_{1,d-1}A_1+\cdots+a_{\rho,d-1}A_{\rho}+(-c_{2,d-1}-c_{1,d-1}c_{2,d})C_2+\cdots+\right.\right.$$
$$\left.\left.(-c_{n,d-1}-c_{1,d-1}c_{n,d})C_n\right) \right)X^{-}_{d-1},g^{-}(A_1)Y^{-}_1,\ldots,g^{-}(A_{\rho})Y^{-}_{\rho},$$
$$\left. g^{-}\left( -(c_{2,d}C_2+\cdots+c_{n,d}C_n)\right) W^{-}_{1},g^{-}(C_2)W^{-}_{2},\cdots,g^{-}(C_n)W^{-}_{n}\right),$$
because $B_1+b_{2,d}B_2+\cdots+b_{m,d}B_m-C_1+c_{2,d}C_2+\cdots+c_{n,d}C_n=0$ on $V$. So by the automorphism $(\hat{x}^{-}_1,\ldots,\hat{x}^{-}_{d-1},\hat{y}^{-}_1,\ldots,\hat{y}^{-}_{\rho},\hat{w})\mapsto(\hat{x}^{+}_1,\ldots,\hat{x}^{+}_{d-1},\hat{y}^{+}_1,\ldots,\hat{y}^{+}_{\rho},\hat{z})$ given by the relations $(11)$ and the action $(12)$, we have $V_t\cong \left(2k,-ka_{1,2},\ldots,-ka_{1,d-1}\right)^{-}V$ for any $t\neq 0$.\hfill q.e.d.

\section{Projective space bundles over the projective line}\label{scroll}

\hspace{.5cm} The classical results for deformations among Hirzebruch surfaces are well-known. As a generalization of this results, for ${\bf P}^{2}$-bundles over ${\bf P}^{1}$, Nakamura \cite{nakamura1} showed the following.

\begin{Prop}[Nakamura \cite{nakamura1}]\label{jisa}

For integers $a,b,c,a',b',c'$, let
$$V={\bf P}_{{\bf P}^1}(\mathcal{O}(a)\oplus\mathcal{O}(b)\oplus\mathcal{O}(c))\mbox{ and }V'={\bf P}_{{\bf P}^1}(\mathcal{O}(a')\oplus\mathcal{O}(b')\oplus\mathcal{O}(c')).$$
Then the following are equivalent.

\begin{enumerate}

\item $a+b+c\equiv a'+b'+c'$ $(${\rm mod} $3)$.
\item There exist ${\bf P}^{2}$-bundles over ${\bf P}^{1}$ $V_0,\ldots,V_m$ such that $V_0\cong V$, $V_m\cong V'$ and $V_{i-1}$ is deformed to $V_{i}$ for any $1\leq i\leq m$.
\item $V$ and $V'$ are homeomorphic.

\end{enumerate}

\end{Prop}

We generalize the case $(1)\Longrightarrow(2)$ of Proposition \ref{jisa} for ${\bf P}^{d-1}$-bundles over ${\bf P}^{1}$ using the one-parameter families constructed in Theorem \ref{maindeform}. Harris \cite{harris1} studied this case. For fundamental properties of primitive collections and primitive relations, see Batyrev \cite{batyrev3}, \cite{batyrev4} and Sato \cite{sato1}. We use the notation as in Section \ref{main}.

Let $V$ be a ${\bf P}^{d-1}$-bundle over ${\bf P}^{1}$, that is,
$$V=V(p_1,\ldots,p_{d-1}):={\bf P}_{{\bf P}^{1}}\left( \mathcal{O}\oplus\mathcal{O}(p_1)\oplus\cdots\oplus\mathcal{O}(p_{d-1})\right),$$
where $p_1,\ldots,p_{d-1}$ are nonnegative integers. Then the primitive relations of the corresponding fan $\Sigma$ are
$${\bf e}_1 + \cdots + {\bf e}_{d-1} + {\bf a}_1=0\ \mbox{and}\ {\bf b}_1+{\bf c}_1=p_1{\bf e}_1+\cdots+p_{d-1}{\bf e}_{d-1},$$
where $\G(\Sigma)=\left\{ {\bf e}_1,\ldots,{\bf e}_{d-1},{\bf a}_1,{\bf b}_1,{\bf c}_1\right\}$. For a nonnegative integer $k$ such that $a_1-k\geq 0$, the conditions in Theorem \ref{maindeform2} are satisfied. Therefore, there exists a one-parameter complex analytic family $\left\{ V_{t}\right\}_{t\in{\bf C}}$ such that
$$(1)\ V_{t}\cong\left\{\begin{array}{ccl}
V & \mbox{if} & t=0 \\
\left(2k,k,\ldots,k\right)^{-}V & \mbox{if} & t\neq 0. \\
\end{array} \right.$$
We show that for $V(p_1,\ldots,p_{d-1})$ and $V(p_1',\ldots,p_{d-1}')$, if $p_1+\cdots+p_{d-1}\equiv p_1'+\cdots+p_{d-1}'\ ($mod$\ d)$, then there exist nonsingular toric $d$-folds $V_0,\ldots,V_m$ such that each $V_i$ is a ${\bf P}^{d-1}$-bundle over ${\bf P}^{1}$, $V_0\cong V(p_1,\ldots,p_{d-1})$, $V_m\cong V(p_1',\ldots,p_{d-1}')$ and $V_{i-1}$ is deformed by a one-parameter family to $V_{i}$ for any $1\leq i\leq m$.

Let $k=1$. Suppose that there exists $1\leq i\leq d-1$ such that $p_i\geq 2$. So we may assume $p_1\geq p_2\geq \cdots \geq p_l>p_{l+1}=\cdots=p_{d-1}=0$ by changing the order of the indices, where $p_1\geq 2$ . Then by the family $(1)$, $V$ is deformed to $(2,1,\ldots,1)^{-}V$. The primitive relations of $(2,1,\ldots,1)^{-}\Sigma$ are
$${\bf e}_1 + \cdots + {\bf e}_{d-1} + {\bf a}_1=0\ \mbox{and}$$
\[
{\bf b}_1+{\bf c}'_1=
\pmatrix{
p_1-2 \cr
p_2-1 \cr
\vdots \cr
p_l-1 \cr
p_{l+1}-1 \cr
\vdots \cr
p_{d-1}-1 \cr
0 \cr}
=
\left\{\begin{array}{ccl}
(p_1-1){\bf e}_1+p_2{\bf e}_2+\cdots+p_{l}{\bf e}_{l}+{\bf a}_1 & \mbox{if} & l<d-1 \\
(p_1-2){\bf e}_1+(p_2-1){\bf e}_2+\cdots+(p_{l}-1){\bf e}_{l} & \mbox{if} & l=d-1, \\
\end{array} \right.$$
where $\G\left( (2,1,\ldots,1)^{-}\Sigma\right)=\left\{ {\bf e}_1,\ldots,{\bf e}_{d-1},{\bf a}_1,{\bf b}_1,{\bf c}_1'\right\}.$ We can replace $V$ by $(2,1,\ldots,1)^{-}V$ and carry out this operation again. This operation terminates in finite steps and $V$ becomes $V(p_1,\ldots,p_{d-1})$ such that $p_1\leq 1,\ldots,p_{d-1}\leq 1$. In each step, $p_1+\cdots+p_{d-1}\in {\bf Z}/d{\bf Z}$ does not change. Thus, we have the following.

\begin{Prop}
For integers $a_1,\ldots,a_d,a'_1,\ldots,a'_d$, let $V={\bf P}_{{\bf P}^1}(\mathcal{O}(a_1)\oplus\cdots\oplus\mathcal{O}(a_d))$ and $V'={\bf P}_{{\bf P}^1}(\mathcal{O}(a'_1)\oplus\cdots\oplus\mathcal{O}(a'_d))$. If $a_1+\cdots+a_d\equiv a'_1+\cdots+a'_d$ $(${\rm mod} $d)$, then there exist ${\bf P}^{d-1}$-bundles over ${\bf P}^{1}$ $V_0,\ldots,V_m$ such that $V_0\cong V$, $V_m\cong V'$ and $V_{i-1}$ is deformed to $V_{i}$ for any $1\leq i\leq m$. Especially, $V$ and $V'$ are homeomorphic.

\end{Prop}

\section{Weakened Fano varieties}\label{weakened}

\hspace{.5cm} The following definition is important for the birational geometry.

\begin{Def}

{\rm
Let $V$ be a nonsingular projective variety. $V$ is called a {\em Fano} (resp. {\em weak Fano}) variety if its anti-canonical divisor $-K_V$ is ample (resp. nef and big).
}

\end{Def}




The following definition was proposed by Minagawa \cite{minagawa1} in connection with ``Reid's fantasy'' for weak Fano $3$-folds.

\begin{Def}[Minagawa \cite{minagawa1}]\label{mina}
{\rm
Let $V$ be a nonsingular weak Fano $d$-fold over ${\bf C}$ and $\Delta_{\epsilon}:=\left\{t\in {\bf C} \, |\, |t|<\epsilon \right\}$ for a sufficiently small real number $\epsilon>0$. Then $V$ is called a {\em weakened Fano} $d$-fold if $V$ is not a nonsingular Fano $d$-fold and there exists a small deformation $\varphi:\mathcal{V}\rightarrow\Delta_{\epsilon}$ such that $\mathcal{V}_0:=\varphi^{-1}(0)\cong V$, while $\mathcal{V}_t:=\varphi^{-1}(t)$ is a nonsingular Fano $d$-fold for any $t\in\Delta_{\epsilon}\setminus\{ 0\}$.
}

\end{Def}

In this section, we give a deformation family for a certain toric weakened Fano $3$-fold using the families constructed in Section \ref{main}. Toric weakened Fano $3$-folds are completely classified by Sato \cite{sato2}. Moreover, we give some examples of toric weakened Fano $4$-folds. We use the notation as in Section \ref{main}.

\begin{Ex}\label{wtofex}
{\rm
Let $V$ be the nonsingular toric weakened Fano $3$-fold of type $X^3_0$ in the sense of Sato \cite{sato2}, that is, the primitive relations of $\Sigma$ are
$${\bf e}_1+{\bf a}_1={\bf e}_2,\ {\bf e}_2+{\bf a}_2=0\mbox{ and }{\bf b}_1+{\bf c}_1=2{\bf e}_1,$$
where $\G(\Sigma)=\left\{ {\bf e}_1,{\bf e}_2,{\bf a}_1,{\bf a}_2,{\bf b}_1,{\bf c}_1\right\}$. $V$ is a $F_1$-bundle over ${\bf P}^1$, where $F_1$ is the Hirzebruch surface of degree $1$. Therefore, by Theorems \ref{maindeform} and \ref{maindeform}, there exists a complex analytic family $\{V_t\}_{t\in{\bf C}}$ such that
\[
\left\{\begin{array}{c}
V_0\cong V,\mbox{ while} \\
V_t\cong (2,-1)^{-}V\ (t\neq 0).\\
\end{array} \right.
\]
The primitive relations of $(2,-1)^{-}\Sigma$ are
$${\bf e}_1+{\bf a}_1={\bf e}_2,\ {\bf e}_2+{\bf a}_2=0\mbox{ and }{\bf b}_1+{\bf c}'_1={\bf e}_2,$$
where $\G(\Sigma)=\left\{ {\bf e}_1,{\bf e}_2,{\bf a}_1,{\bf a}_2,{\bf b}_1,{\bf c}'_1\right\}$. $(2,-1)^{-}V$ is the toric Fano $3$-fold we want (see Section $4$ in Sato \cite{sato2}).
}

\end{Ex}

In the same way as in Example \ref{wtofex}, we get some examples of toric weakened Fano $4$-folds which does not decomposed into direct products of lower-dimensional varieties. In the following, $\G(\Sigma)=\{x_1,x_2,\ldots\}$ and the fans corresponding to toric weakened Fano $4$-folds are described in terms of primitive relations. We also give the types of general fibers. The types of nonsingular toric Fano $4$-folds are in the sense of Batyrev \cite{batyrev4} and Sato \cite{sato1}.


\begin{enumerate}
\item $x_1+x_4=x_2,\ x_2+x_3+x_5=0\ $and$\ x_6+x_7=2x_1$ (type $D_7$).
\item $x_1+x_4=x_2,\ x_2+x_6=0,\ x_3+x_5=x_2\ $and$\ x_7+x_8=2x_1$ (type $L_1$).
\item $x_1+x_4=x_2,\ x_2+x_6=0,\ x_3+x_5=x_6\ $and$\ x_7+x_8=2x_1$ (type $L_{13}$).
\item $x_1+x_4=x_2,\ x_2+x_5=x_3,\ x_3+x_6=0\ $and$\ x_7+x_8=2x_1$ (type $L_2$).
\item $x_5+x_6=0,\ x_3+x_7=0,\ x_2+x_3=x_5,\ x_5+x_7=x_2,\ x_2+x_6=x_7,\ x_1+x_4=x_2\ $and$\ x_8+x_9=2x_1$ (type $Q_1$).
\item $x_5+x_6=0,\ x_3+x_7=0,\ x_2+x_3=x_5,\ x_5+x_7=x_2,\ x_2+x_6=x_7,\ x_1+x_4=x_3\ $and$\ x_8+x_9=2x_1$ (type $Q_{13}$).
\item $x_5+x_6=0,\ x_3+x_7=0,\ x_2+x_3=x_5,\ x_5+x_7=x_2,\ x_2+x_6=x_7,\ x_1+x_4=x_5\ $and$\ x_8+x_9=2x_1$ (type $Q_8$).
\item $x_5+x_6=0,\ x_3+x_7=0,\ x_2+x_3=x_5,\ x_5+x_7=x_2,\ x_2+x_6=x_7,\ x_1+x_4=0\ $and$\ x_8+x_9=2x_1$ (type $Q_{11}$).
\item $x_5+x_8=0,\ x_2+x_5=x_3,\ x_3+x_8=x_2,\ x_3+x_6=x_5,\ x_3+x_7=0,\ x_2+x_6=0,\ x_6+x_8=x_7,\ x_2+x_7=x_8,\ x_5+x_7=x_6,\ x_1+x_4=x_2\ $and$\ x_9+x_{10}=2x_1$ (type $U_1$).
\end{enumerate}

\bigskip

\begin{flushleft}
\begin{sc}
Department of Mathematics,\\
Tokyo Institute of Technology,\\
2-12-1 Oh-okayama, Meguro-ku, Tokyo 152-8551, Japan.
\end{sc}

\medskip
{\it E-mail address}: $\mathtt{hirosato@math.titech.ac.jp}$
\end{flushleft}

\end{document}